# A Simple Method for Generating Rational Triangles

*Konstantine Zelator*
*Department Of Mathematics*
*College Of Arts And Sciences*
*Mail Stop 942*
*University Of Toledo*
*Toledo,OH 43606-3390*
*U.S.A.*

# 1 Introduction

The purpose of this paper is to present a simple method for generating rational triangles. These are triangles with rational side lengths and rational area. We do this in Section 5, wherein we generate a certain family of such triangles. In Section 2, we present some well known triangle formulas involving a triangle's side lengths, heights, area, the Law of Sines, and the Law of Cosines. Also two formulas for $R$ and $r$, the radius of a triangle's circumscribed circle and the radius of its inscribed circle respectively. We offer a short derivation of these two formulas in Section 7. In Section 3, we state the well known parametric formulas for Pythagorean triples, as well as a table featuring a few Pythagorean triples. In Section 4 we state two definitions; that of a rational triangle and the definition of a Pythagorean rational. And in Section 6, we offer two closing remarks.

# 2 Triangle Formulas

For a triangle $AB\mathrm{G}$, let $a = |\overline{B\mathrm{G}}|, b = |\overline{A\mathrm{G}}|, g = |\overline{B\mathrm{G}}|$, be its three side lengths; $h_a, h_b, h_g$ the three heights (corresponding to the perpendiculars from the vertex $A$ to the side $\overline{B\mathrm{G}}$, from $B$ to $\overline{A\mathrm{G}}$, and from G to $\overline{AB}$ respectively); $E$ the triangle's area. Also let $R$ be the radius of the triangle's circumscribed circle, and $r$ the radius of the triangle's inscribed circle. The following four sets of formulas are very well known to a wide mathematical audience.

(2)

$$E = \frac{1}{2} bg \sin A = \frac{1}{2} ag \sin B = \frac{1}{2} ab \sin \mathrm{G} \qquad (1)$$

(3)

Law of sines: $\frac{a}{\sin A} = \frac{b}{\sin B} = \frac{g}{\sin G} = 2R$

Law of cosines: $\cos A = \frac{b^2 + g^2 - a^2}{2bg}, \cos B = \frac{a^2 + g^2 - b^2}{2ag}, \cos G = \frac{a^2 + b^2 - g^2}{2ab}$

(4)

A notational remark here: We have used the single Greek letters $A, B, G$ to denote the triangle $ABG$'s internal angles (in other words $A$ for instance, stands for the angle $BAG$; or similarly for the angle $GAB$. Of course, the same single capital letters also denote the triangles vertices.

The next two formulas are not as well known to a wide mathematical audience as the first four.

(5)

$$R = \frac{abg}{4E}$$

(6)

$$r = \frac{2E}{a + b + g}$$

We present a short derivation for these two formulas in Section 7.

# 3 Pythagorean Triples

A Pythagorean triple *(a, b,* c) of positive integers *a, b, c*; is one which satisfies

$a^2 + b^2 = c^2$. In effect, the integers *a* and b are the leg lengths of a right triangle;

while $c$ is its hypotenuse length. Such a right triangle is widely known as a Pythagorean triangle. If ($a, b, c$) is a Pythagorean triple, then

$a = d(m^2 - n^2), b = d.(2mn), c = d.(m^2 + n^2)$, for some positive integers $d, m, n$ such that $m > n$, $(m,n) = 1$(i.e., $m$ and $n$ are relatively prime, and $m + n \equiv 1 \pmod 2$ (that is, $m$ and $n$ have different parities; one is odd, the other even). (Alternatively, of course, we can have $a = d(2mn)$, and $b = d(m^2 - n^2)$, instead). (7)

Conversely, any triple ($a, b, c$) which satisfies (7) must be a Pythagorean triple, as a straightforward calculation shows. In effect, formulas (7) describe the entire family of Pythagorean triples.

When $d = 1$ the resulting Pythagorean triple is called *primitive.*

This material can be found in almost any number theory text. For example, refer to [1] for such a text. L.E. Dicksons monumental book (see [2]) has a wealth of historical information on the subject. Below, we present a listing of the first six primitive Pythagorean triples: those obtained from formulas (7) with $d = 1$ and under the constraint $m \le 5$.

| $d$ | m | n | a | b | c |
|---|---|---|---|---|---|
| 1 | 2 | 1 | 3 | 4 | 5 |
| 1 | 3 | 2 | 5 | 12 | 13 |
| 1 | 4 | 1 | 15 | 8 | 17 |
| 1 | 4 | 3 | 7 | 24 | 25 |
| 1 | 5 | 2 | 21 | 20 | 29 |
| 1 | 5 | 4 | 9 | 40 | 41 |

**Table 1**

# 4 Two Definitions

**Definition 1:** A rational number $r$ is called a *Pythagorean rational,* if $r = \frac{a}{b}$, where $a$ and $b$ are relatively prime positive integers such that $a^2 + b^2 = c^2$, for some positive integer $c$.

**Remark 1:** If r is a Pythagorean rational, then obviously, so is the rational number $\frac{1}{r}$; since $a^2 + b^2 = c^2 = b^2 + a^2$.

**Remark 2:** Note that the positive integers $a$ and $b$ in Definition 1, are the leg lengths of a primitive Pythagorean triangle whose hypotenuse length is $c$; $(a, b, c)$ is a primitive Pythagorean triple.

**Definition 2;** A triangle $AB$G with side lengths $a, b, g$ and area $E$, is said to be a *rational triangle,* if all four $a, b, g$ and $E$, are rational numbers.

**Remark 3:** It is easy to come up with a triangle whose side lengths are rational numbers. Indeed, all one has to do is pick three positive rationals $a, b, g$, that satisfy the three triangle inequalities $a + b > g, b + g > a$ , and $a + g > b$. To ensure that the triangle is non-isosceles, the three rationals must be distinct. And if we want it to be non-right, we must ensure that $a^2 + b^2 \pi g^2, a^2 + g^2 \pi b^2$, and $a^2 \pi b^2 + g^2$ . On the other hand, to find by trial and error, a non-right non-isosceles rational triangle is significantly more difficult. One would have to find positive rationals, which in addition to the above conditions, must be such that the product $(a + b + g)(-a + b + g)(a - b + g)(a + b - g)$ is the square of a rational number. This is because the area $E$ satisfies the formula,

$$E = \frac{1}{4}.\sqrt{(\alpha+\beta+\gamma)(-\alpha+\beta+\gamma)(\alpha-\beta+\gamma)(\alpha+\beta-\gamma)},$$

well-known in the literature as *Herons formula.*

*Observations*

By inspection of formuals (1) - (6) the following is evident.

1. If $\alpha, \beta, \gamma$ are rational, then the three cosine values $\cos A, \cos B, \cos \Gamma$, are always rational, regardless of whether the triangle $AB\Gamma$ is rational.

2. If $AB\Gamma$ is a rational triangle, then the three heights, the three sine values, and the radii $R$ and $\rho$ are also rational.

3. If $AB\Gamma$ is not a rational triangle, then $E$ is an irrational number, and consequently the three sine values, the three heights, and the two radii are all irrational numbers.

**Remark 4:** Note that a Pythagorean rational as defined in Definition 1, is always a *proper rational;* that is, it cannot be an integer. This can be easily seen from the conditions and formulas in (7). On the other hand, the notion of a *Pythagorean number,* is often found in the literature of number theory.

A Pythagorean number is a number which is equal to the area of a Pythagorean triangle. Thus, by (7), it is always an integer equal to $\delta^2.mn.(m^2 - n^2)$ a Pythagorean number is never a proper rational.

# 5 A Family of Rational Triangles

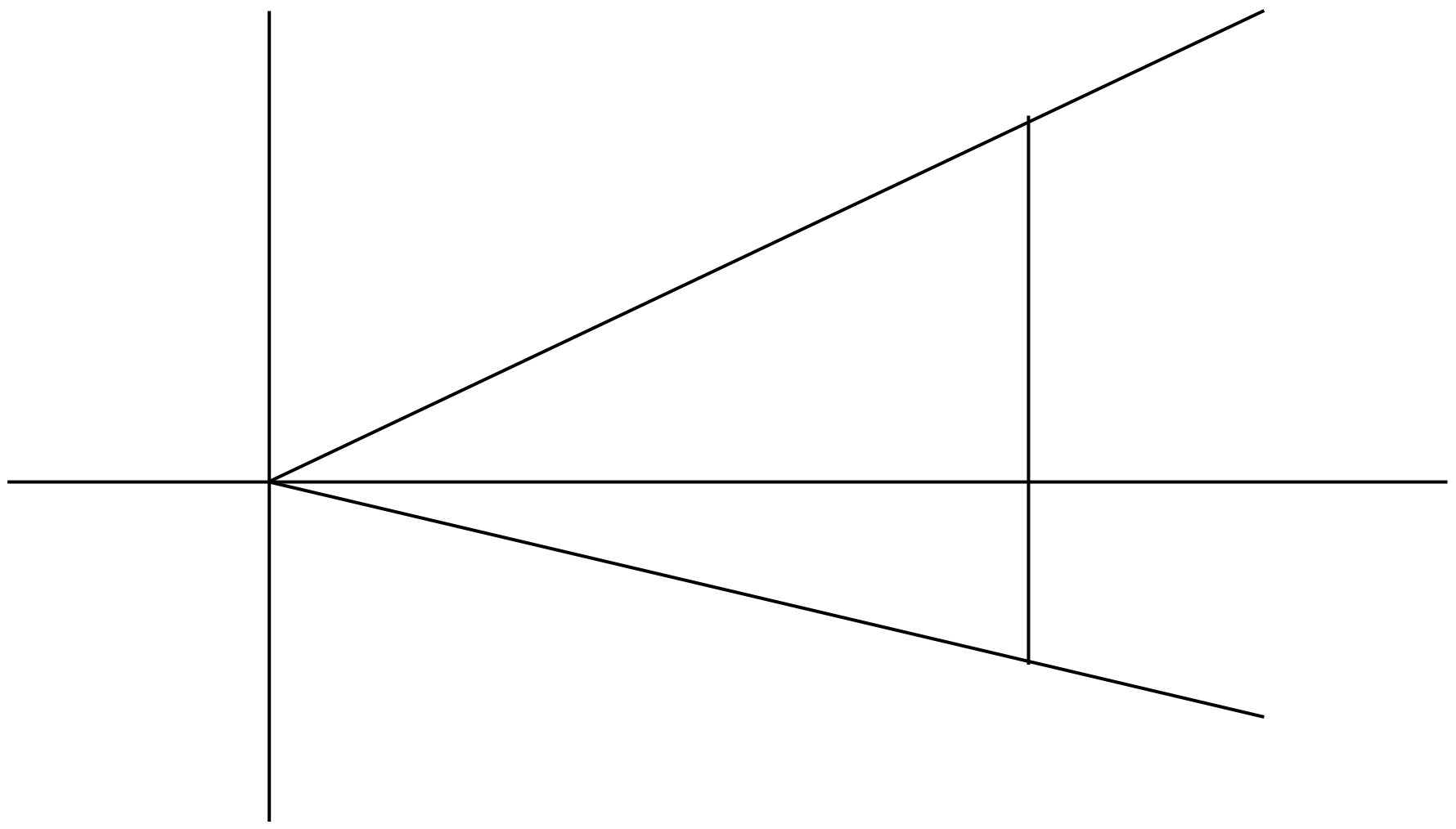

**Figure 1**

Let $r_1 = \frac{a_1}{b_1}$ and $r_2 = \frac{a_2}{b_2}$ be two Pythagorean rationals, which means that $a_1^2 + b_1^2 = c_1^2$, and $a_2^2 + b_2^2 = c_2^2$ for some positive integer $c_1$ and $c_2$.

In the $x$- $y$ plane, consider two rays emanating from the origin O(0,0). One ray lies in the first quadrant and has slope $m_1 = \frac{a_1}{b_1}$, the other ray lies in the fourth quadrant and has slope $m_2 = -\frac{a_2}{b_2}$. Also consider the points $A_1(1, m_1)$ and $A_2(1, m_2)$ on the first and second rays respectively; T(l, 0) on the positive o>axis.

The triangle $A_1OA_2$ is a rational triangle. We put,

$\varphi_1$ = degree measure of acute angle $\angle A_1OT$; $0 < \varphi_1 < 90$

$\varphi_2$ = degree measure of acute angle $\angle A_2OT$; $0 < \varphi_2 < 90$

$\varphi = \varphi_1 + \varphi_2$ = degree measure of angle $\angle A_1OA_2$; $0 < \varphi < 180$

$\omega_1$ = degree measure of acute angle $\angle A_1OA_2$; $0 < \omega_1 < 90$

$\omega_2$ = degree measure of acute angle $\angle A_2OA_1$; $0 < \omega_2 < 90$

Also for side lengths: $\alpha = \left|\overline{OA_1}\right|, \beta = \left|\overline{OA_2}\right|, \gamma = \left|\overline{A_1A_2}\right|$.

Heights: $h_\alpha, h_\beta, h_\gamma$.

(8)

Straightforward calculations yield the following results.

*Side lengths*

$$\alpha = \sqrt{1+m_1^2} = \frac{c_1}{b_1}, \beta = \sqrt{1+m_2^2} = \frac{c_2}{b_2}, \gamma = m_1 - m_2 = \frac{a_1b_2 + a_2b_1}{b_1b_2}$$

*Area*

$$\text{Area } E = \frac{1}{2}.1.(m_1 - m_2)^2 = \frac{a_1b_2 + a_2b_1}{2b_1b_2} \qquad (9)$$

*Heights*

$$h_\alpha = \frac{2E}{\alpha} = \frac{a_1b_2 + a_2b_1}{c_1b_2}, h_\beta = \frac{2E}{\beta} = \frac{a_1b_2 + a_2b_1}{c_2b_1}, h_\gamma = 1 \qquad (10)$$

*Cosine and Sine Values*

$$
\begin{aligned}
&\cos\varphi_1 = \frac{b_1}{c_1}, \sin\varphi_1 = \frac{a_1}{c_1}, \cos\varphi_2 = \frac{b_2}{c_2}, \sin\varphi_2 = \frac{a_2}{c_2} \\
&\cos\varphi = \cos(\varphi_1 + \varphi_2) = \cos\varphi_1 \cos\varphi_2 - \sin\varphi_1 \sin\varphi_2 = \frac{b_1 b_2 - a_1 a_2}{c_1 c_2} \\
&\sin\varphi = \sin(\varphi_1 + \varphi_2) = \sin\varphi_1 \cos\varphi_2 + \sin\varphi_2 \cos\varphi_1 = \frac{a_1 b_2 + a_2 b_1}{c_1 c_2} \\
&\cos\omega_1 = \frac{a_1}{c_1}, \sin\omega_1 = \frac{b_1}{c_1}, \cos\omega_2 = \frac{a_2}{c_2}, \sin\omega_2 = \frac{b_2}{c_2}
\end{aligned}
\tag{11}
$$

$R$ and $\rho$

(To obtain these we use (5) and (6), combined with (8) and (9))

$$R = \frac{c_1 c_2}{2 b_1 b_2}, \rho = \frac{2(a_1 b_2 + a_2 b_1)}{c_1 b_2 + b_1 c_2 + a_1 b_2 + a_2 b_1}$$

*Observations*

4. One can easily see from Figure 1, that the triangle $OA_1A_2$ will be a non-right one, if and only if $a_1a_2 \neq b_2b_1$. Equivalently, the same triangle will be a right triangle, precisely when $a_1a_2 = b_2b_1$. This couched with the coprimeness conditions $(a_1b_1) = 1 = (a_2, b_2)$, and a bit of number theory, easily implies that $a_1 = b_2$ and $b_1 = a_2$. In other words, $OA_1A_2$ will be a right triangle, exactly when the two Pythagorean rationals are reciprocals of each other (and so, $m_1m_2 = -1$).

5. Also easily, by inspection from Figure 1, one sees that the triangle $A_1OA_2$ will be non-isosceles if and only if, $a_1b_2 \neq a_2b_1$. Equivalently, it will be isosceles precisely when $a_1b_2 = a_2b_1$ which when combined with the coprimeness conditions, it easily implies $a_1 = a_2$ and $b_1 = b_2$; the two Pythagorean rationals $r_1$ and $r_2$ are equal in this case.

If two distinct Pythagorean rationals are given, say $r$ and s; then eight different triangles $A_1OA_2$ can be formed as follows:

Triangle 1: By taking $m_1 = r$ and $m_2 = -s$

Triangle 2: By taking $m_1 = r$ and $m_2 = -\frac{1}{s}$

Triangle 3: By taking $m_1 = s$ and $m_2 = -r$

Triangle 4: By taking $m_1 = s$ and $m_2 = -\frac{1}{r}$

Triangle 5: By taking $m_1 = \frac{1}{r}$ and $m_2 = -s$

Triangle 6: By taking $m_1 = \frac{1}{r}$ and $m_2 = -\frac{1}{s}$

Triangle 7: By taking $m_1 = \frac{1}{s}$ and $m_2 = -r$

Triangle 8: By taking $m_1 = \frac{1}{s}$ and $m_2 = -\frac{1}{r}$

However, Triangle 3 is simply the perpendicular reflection of Triangle 1 through the x-axis, as it is evident from Figure 1. Thus, Triangles 1 and 3 are congruent. Likewise Triangles 2 and 7 are congruent. Triangles 4 and 5, and Triangles 6 and 8 are also congruent. In summary, the group of eight triangles, can be divided into four sets, each set containing a pair of congruent triangles, because of reflection through the x-axis. Now, let us take a look at Table 1, which contains six rows of numbers. If we fix two of the six rows; we construct eight triangles $A_1OA_2$ as the above analysis shows. In how many ways can we choose two rows from the given six? Every combinatorics student knows that the answer to this is given by the binomial coefficient,

$$\binom{6}{2} = \frac{6!}{2!4!} = \frac{5.6}{2} = 15.$$

Thus, by implementing the above method, we see that there are 120 such triangles that arise from Table 1; sixty pairs of triangles; each pair containing two congruent triangles.

# 6 A Closing Remark

A glance at formulas (7) with $d = 1$ shows that,

$$\frac{a}{b} = \frac{1}{2}\left(\frac{m}{n} - \frac{n}{m}\right), \quad \frac{b}{a} = \frac{2mn}{m^2 - n^2}$$

Thus, if we take n = 1 (and *m* even), it is obvious that the Pythagorean rationals $\frac{a}{b}$ can be arbitrarily large; while $\frac{b}{a}$ arbitrarily small. Therefore in Figure 1, either $m_1$ or $m_2$ can be arbitrarily large or small; and so the degree measure $J$ of the angle $A_1OA_2$ can be arbitrarily close to 180 or to 0.

# 7 A Derivation of Formulas (5) and (6)

The derivation of (5) is very short indeed. From (2) and (3) we obtain

$$E = \frac{1}{2}.bg \sin A; \quad E = \frac{1}{2} bg.\frac{a}{2R} \;\hat{U}\; R = \frac{abg}{4E}$$

The derivation of (6) is also very short, if we just take a look at Figure 2. We have

$$\text{Area}\,(ABG) = E = \text{Area}\,(AIB) + \text{Area}\,(BIG) + \text{Area}\,(GIA);$$

$$E = \frac{1}{2}.r.g + \frac{1}{2}.r.a + \frac{1}{2}.r.b;$$

$$E = \frac{1}{2}.r(g + a + b); \; r = \frac{2E}{a + b + g}$$

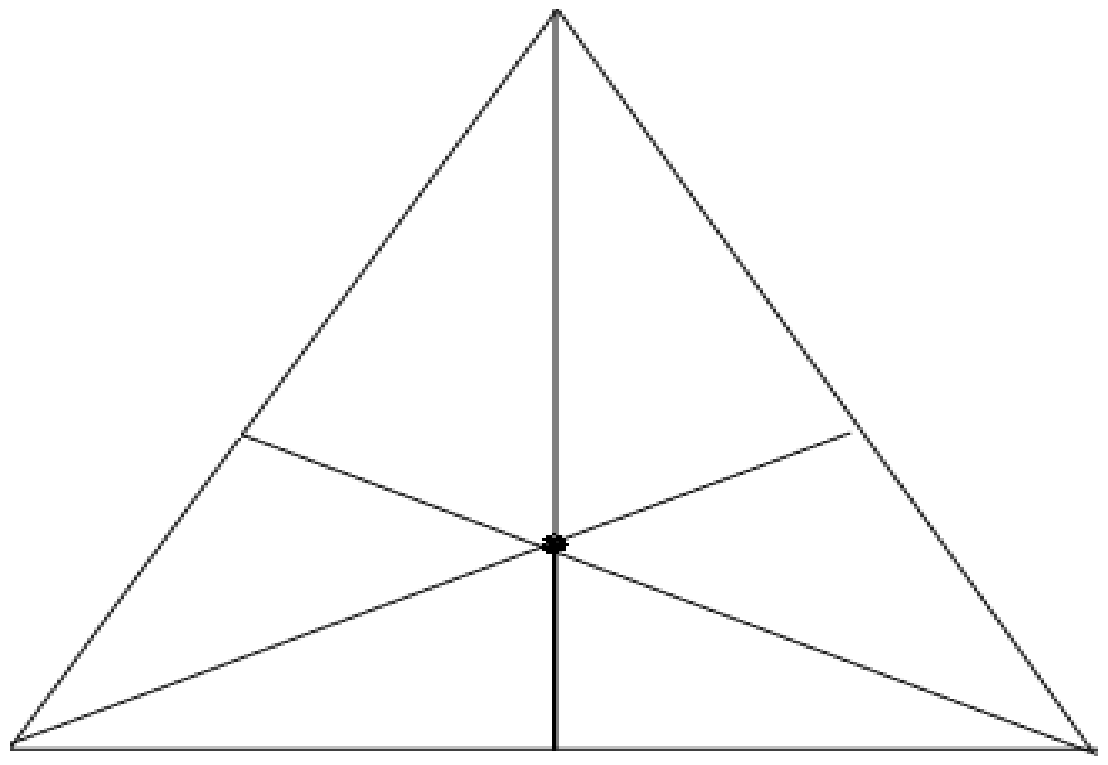

**Figure 2**

# 8 References

1. Rosen, Kenneth H., *Elementary Number Theory and its Applications,* Third edition, 1993, Addison-Wesley Publishing Company (there is now a fourth edition as well), 544 pp. ISBN:0-201-57889-l, For Pythagorean triples, see pages 436-442.

2. Dickson, L.E., *History of Theory of Numbers,* Vol. II, AMS Chelsea Publishing, Providence, Rhode Island, 1992. ISBN: 0-8218-1935-6; 803 pp (unaltered textual reprint of the original book, first published by Carnegie Institute of Washington in 1919, 1920, and 1923). (a) For material on Pythagorean triangles and rational right triangles, see pages 165-190. (b) For material on rational triangles in general, see pages 191-216.